\documentclass[a4paper]{article}
\usepackage{geometry}
\usepackage{amssymb}
\usepackage{latexsym}
\usepackage{amsmath}
\usepackage{indentfirst}
\usepackage{graphicx}
\usepackage[colorlinks=true]{hyperref}
\usepackage{mathrsfs} 
\usepackage{chngcntr}
\usepackage{color}

\newtheorem{theorem}{Theorem}[section]

\newtheorem{lemma}{Lemma}[section]

\numberwithin{equation}{section}

\newenvironment{proof}{\medskip\par\noindent{\bf Proof\/}:\quad}{\qquad
	\raisebox{-0.5mm}{\rule{1.5mm}{1mm}}\vspace{6pt}}

\title{Existence of normalized solutions to a class of fractional Choquard equations with potentials}

\author{Yongpeng Chen, Zhipeng Yang\thanks{Corresponding author: yangzhipeng326@163.com}, Jianjun Zhang}    

\date{}

\AtEndDocument{%
	\par
	\bigskip
	\bigskip
	
	\noindent
	\textbf{Yongpeng Chen}\\[0.2em]
	\textsc{Yunnan Key Laboratory of Modern Analytical Mathematics and Applications, Yunnan Normal University, Kunming, China}\\[0.3em]
	\textsc{School of Science, Guangxi University of Science and Technology, Liuzhou, China.}\\[0.3em]
	\textit{E-mail address}: \texttt{ yongpengchen@mail.bnu.edu.cn}\\[1em]
	
	\noindent
	\textbf{Zhipeng Yang}\\[0.2em]
	\textsc{Department of Mathematics, Yunnan Normal University, Kunming, China}\\
	\textsc{Yunnan Key Laboratory of Modern Analytical Mathematics and Applications, Yunnan Normal University, Kunming, China}\\[0.3em]
	\textit{E-mail address}: \texttt{yangzhipeng326@163.com}
	
	\bigskip
	
	\noindent
	\textbf{Jianjun Zhang}\\[0.2em]
	\textsc{College of Mathematics and Statistics, Chongqing Jiaotong University, Chongqing, China}\\[0.3em]
	\textit{E-mail address}: \texttt{zhangjianjun09@tsinghua.org.cn}\\[1.5em]
}

\begin{document}
	\date{} \maketitle
	\begin{abstract}
This paper investigates the existence of normalized solutions to the nonlinear fractional Choquard equation:
$$
(-\Delta)^s u+V(x) u=\lambda u+f(x)\left(I_\alpha *\left(f|u|^q\right)\right)|u|^{q-2} u+g(x)\left(I_\alpha *\left(g|u|^p\right)\right)|u|^{p-2} u, \quad x \in \mathbb{R}^N
$$
subject to the mass constraint
$$
\int_{\mathbb{R}^N}|u|^2 d x=a>0,
$$
where $N>2 s, s \in(0,1), \alpha \in(0, N)$, and $\frac{N+\alpha}{N} \leq q<p \leq \frac{N+\alpha+2 s}{N}$. Here, the parameter $\lambda \in \mathbb{R}$ appears as an unknown Lagrange multiplier associated with the normalization condition. By employing variational methods under appropriate assumptions on the potentials $V(x), f(x)$, and $g(x)$, we establish several existence results for normalized solutions.
	\end{abstract}

\section{Introduction}
In this paper, we investigate the existence of normalized solutions to the following nonlinear fractional Choquard equation:
\begin{equation}\label{eq1}
	\begin{cases}
		(-\Delta)^s u + V(x)u = \lambda u + f(x)\left(I_\alpha *(f|u|^q)\right)|u|^{q-2} u + g(x)\left(I_\alpha *(g|u|^p)\right)|u|^{p-2} u,\\
		\int_{\mathbb{R}^N}|u|^2 d x=a,
	\end{cases}
\end{equation}
where $N>2 s, s \in(0,1), \alpha \in(0, N)$ and
$$
\frac{N+\alpha}{N} \leq q<p \leq \frac{N+2 s+\alpha}{N}.
$$
The parameter $\lambda \in \mathbb{R}$ arises as an unknown Lagrange multiplier associated with the mass constraint. The nonlocal interaction terms involve the Riesz potential $I_{\alpha}$ which is defined by
$$
I_\alpha(x)=\frac{A_\alpha}{|x|^{N-\alpha}}, \quad x \in \mathbb{R}^N \backslash\{0\}, \quad \text { with } \quad A_\alpha=\frac{\Gamma\left(\frac{N-\alpha}{2}\right)}{\Gamma\left(\frac{\alpha}{2}\right) \pi^{\frac{N}{2}} 2^\alpha},
$$
where $\Gamma$ denotes the Gamma function. The fractional Laplace operator $(-\Delta)^s$ is given by the principal value integral
$$
(-\Delta)^s u(x)=C(N, s) \text { P.V. } \int_{\mathbb{R}^N} \frac{u(x)-u(y)}{|x-y|^{N+2 s}} d y, \quad x \in \mathbb{R}^N.
$$
where $C(N, s)$ is a normalization constant (see \cite{MR2944369}). Introduced by Laskin \cite{MR1755089} in the context of Feynman path integrals associated with L\'evy-type quantum processes, the fractional Laplacian plays an essential role in modeling phenomena such as Boson stars and water-wave dynamics \cite{MR2318846}.
\par 
The analysis of normalized solutions to the fractional Choquard equation
$$
(-\Delta)^s u+V(x)u=\lambda u+\left(I_\alpha *|u|^p\right)|u|^{p-2} u, \quad x \in \ \mathbb{R}^N,
$$
depends fundamentally on three critical exponents:
\begin{itemize}
	\item the Hardy-Littlewood-Sobolev lower critical exponent $$\frac{N+\alpha}{N};$$
	\item the $L^2$-critical exponent $$\frac{N+2 s+\alpha}{N};$$
	\item the Hardy-Littlewood-Sobolev upper critical exponent $$2_{\alpha, s}^*=\frac{N+\alpha}{N-2 s}.$$
\end{itemize} 
The location of the exponent $p$ with respect to these critical thresholds determines the appropriate analytical framework and shapes the technical challenges encountered in each scenario.
\par 
D'Avenia et al. \cite{Avenia} studied 
\[
(-\Delta)^s u+\omega u=\left(I_\alpha *|u|^p\right)|u|^{p-2} u, \quad u \in H^s\left(\mathbb{R}^N\right),
\]
on the constraint
$$
S(c)=\left\{u \in H^s\left(\mathbb{R}^N\right):\int_{\mathbb{R}^N}|u|^2 d x=c\right\},
$$
where $N \geqslant 3, s \in(0,1), \alpha \in(0, N), p \in\left(1+\frac{\alpha}{N}, 1+\frac{\alpha+2 s}{N}\right)$. The associated energy functional is
\begin{equation*}
	E(u)=\frac{1}{2}\int_{\mathbb{R}^N}|(-\Delta)^{\frac{s}{2}} u|^2dx-\frac{1}{2p} \int_{\mathbb{R}^N}\left(I_\alpha *|u|^p\right)|u|^p d x.
\end{equation*}
For any \( c > 0 \), the existence of a normalized solution was established by considering a minimizing sequence of the functional \( E \) over the constraint set \( S(c) \).  Subsequently, Li and Luo \cite{li-luo} extended their existential results to the mass-supercritical regime. Specifically, for 
$
p \in \left( \max\left\{1 + \frac{\alpha + 2s}{N}, 2\right\}, \frac{N + \alpha}{N - 2s} \right),
$
they employed a constrained minimization method on a Pohozaev-Nehari manifold associated with \( S(c) \). Their analysis proved that for every \( c > 0 \), the problem admits a normalized
ground state solution. Additionally, they investigated the asymptotic behavior of these solutions in the limits \( c \to 0^+ \) and \( c \to +\infty \).  Beyond ground state solutions, Li and Luo also utilized a minimax approach to establish the existence of infinitely many high-energy normalized solutions for every \( c > 0 \).

He et al. \cite{he-zou} investigated the following fractional Choquard equation with an upper Hardy-Littlewood-Sobolev  critical exponent:  
\begin{equation*}
	\begin{cases}(-\Delta)^s u=\lambda u+\mu|u|^{q-2} u+\left(I_\alpha *|u|^{2_{\alpha, s}^*}\right)|u|^{2_{\alpha, s}^*-2} u, \quad x \in \mathbb{R}^N, \\ 
		\int_{\mathbb{R}^N}|u|^2 d x=a^2,
	\end{cases}
\end{equation*}
where  $s \in(0,1), N>$ $2 s, 0<\alpha<\min \{N, 4 s\}, 2<q<\frac{2 N}{N-2 s}$, and  $ \mu \in \mathbb{R}$. In the \( L^2 \)-subcritical regime, they proved the existence of a ground state solution corresponding to a local minimizer of the energy functional.  For the \( L^2 \)-critical and \( L^2 \)-supercritical cases, ground states were constructed as mountain pass solutions.  
Furthermore, the behavior of ground state solutions was studied as \( \mu \to 0 \).

Chen and Yang \cite{chen-yang} examined the following problem with a lower Hardy-Littlewood-Sobolev  critical exponent: 
$$
\left\{\begin{array}{l}
	(-\Delta)^s u+\lambda u=\gamma\left(I_\alpha *|u|^{\frac{\alpha}{N}+1}\right)|u|^{\frac{\alpha}{N}-1} u+\mu\left(I_\alpha *|u|^q\right)|u|^{q-2} u, \quad \text { in } \mathbb{R}^N, \\
	\int_{\mathbb{R}^N}|u|^2 \mathrm{~d} x=c^2,
\end{array}\right.
$$
where $N \geq 3, s \in(0,1), \alpha \in(0, N), \gamma, \mu, c>0$ and $\frac{N+\alpha}{N}<q \leq 2_{\alpha, s}^*$. The existence and non-existence of normalized solutions to this problem were established, depending on which interval among those defined by $\frac{N+\alpha}{N}$, $\frac{N+\alpha+2s}{N}$, and $2_{\alpha, s}^*$ the parameter $q$ lies in. 

Recently, Chen, Yang and Zhang \cite{ChenYangZhangNoDEA} studied a fractional Choquard equation with mixed nonlinearities and a small parameter,
\[
(-\Delta)^s u+V(\varepsilon x)u
=\lambda u+\left(I_\alpha*|u|^q\right)|u|^{q-2}u
+\left(I_\alpha*|u|^p\right)|u|^{p-2}u
\quad \text{in } \mathbb{R}^N,
\]
under the prescribed mass constraint
\[
\int_{\mathbb{R}^N}|u|^2\,dx=a>0.
\]
Here the exponents satisfy
\[
\frac{N+\alpha}{N}<q<\frac{N+2s+\alpha}{N}<p
\leq \frac{N+\alpha}{N-2s}.
\]
By combining variational methods with the Lusternik--Schnirelmann category theory, they obtained multiple normalized solutions whose number is related to the topology of the set where the potential attains its minimum values. 

Liu et al.\cite{liu-sun-zhang}  studied the following non-autonomous fractional Choquard equation with a potential term:
$$
\left\{\begin{array}{l}
	(-\Delta)^s u-\left(I_\alpha *\left[h(\varepsilon y)|u|^{\frac{N+\alpha}{N}}\right]\right) h(\varepsilon x)|u|^{\frac{N+\alpha}{N}-2} u-\mu|u|^{q-2} u=\lambda u,\  \text { in } \mathbb{R}^N, \\
	\int_{\mathbb{R}^N} u^2 d x=c,
\end{array}\right.
$$
where $s \in(0,1), \alpha \in(0, N), N>2 s, 2<q<2+\frac{4 s}{N}, \mu, c>0, \lambda \in \mathbb{R}$ appears as an unknown Lagrange multiplier.
The function $h$ satisfies:
\begin{itemize}
	\item[$(h_1)$] $h \in C\left(\mathbb{R}^N\right), 0<\inf _{x \in \mathbb{R}^N} h(x)=\lim _{|x| \rightarrow+\infty} h(x) \leq \max _{x \in \mathbb{R}^N} h(x)=: h_{\text {max }} ;$
	\item[$(h_2)$] $h^{-1}\left(h_{\max }\right)=\left\{a_1, a_2, \cdots, a_l\right\}$ with $a_1=0$ and $a_i \neq a_j$ if $i \neq j$.
\end{itemize}
Under hypothesis $(h_1)$, the authors established the existence of at least one ground state solution. Furthermore, when both $(h_1)$ and $(h_2)$ hold, they showed that for $\varepsilon$ sufficiently small, the number of distinct normalized solutions is bounded below by the number of global maximum points of  $h$.

For additional results concerning normalized solutions of the Choquard equation and related fractional Choquard problems, we refer to the works \cite{ChenKumarYangZhang2026,ChenLiYang2020,MR4744303,MR4576453,MR4707368,MR4767698,MR4668062,MR4526928,MR4579263,YangZhao2021,Yu2025NormalizedSF,MR5035604,MR4854878}.
To date, research on this subject has primarily focused on autonomous problems. However, the presence of a potential can fundamentally alter the existence of solutions, a phenomenon well-established in the study of Schr\"{o}dinger equations, yet remains largely unexplored for the fractional Choquard equation.  Motivated by these considerations, we investigate the existence of normalized solutions to Problem \eqref{eq1} using variational methods. By imposing suitable structural conditions on the potentials, we establish rigorous existence theorems.  Before stating our main results, we first introduce some notations

For any $s \in(0,1)$, the Sobolev space $H^s\left(\mathbb{R}^N\right)$ is defined as:
$$
H^s\left(\mathbb{R}^N\right) =\left\{u \in L^2\left(\mathbb{R}^N\right): (-\Delta)^{\frac{s}{2}} u\in L^2\left(\mathbb{R}^N\right)\right\}
$$
with the norm:
$$
\|u\|_{H^s\left(\mathbb{R}^N\right)}=\left(\int_{\mathbb{R}^N}|u|^2dx+\int_{\mathbb{R}^N}|(-\Delta)^{\frac{s}{2}} u|^2dx\right)^{\frac{1}{2}},
$$
and the associated inner product:  
$$
\left\langle u, v\right\rangle=\int_{\mathbb{R}^N}(-\Delta)^{\frac{s}{2}} u(-\Delta)^{\frac{s}{2}} v d x+\int_{\mathbb{R}^N} uvdx .
$$
Additionally, we introduce the space:  
\[\mathcal{D}^{s, 2}\left(\mathbb{R}^N\right)=\left\{u \in L^{2_s^*}\left(\mathbb{R}^N\right): (-\Delta)^{\frac{s}{2}} u\in L^2\left(\mathbb{R}^N\right)\right\},
\]
with the norm 
\[
\|u\|=\left(\int_{\mathbb{R}^N}|(-\Delta)^{\frac{s}{2}} u|^2dx\right)^{\frac{1}{2}}.
\]
For $u \in H^s\left(\mathbb{R}^N\right)$, applying Proposition  3.6 from  \cite{MR2944369}, we obtain:
$$
\int_{\mathbb{R}^N}|(-\Delta)^{\frac{s}{2}} u|^2dx=\frac{1}{2} C(N, s) \int_{\mathbb{R}^N} \int_{\mathbb{R}^N}\frac{|u(x)-u(y)|^2}{|x-y|^{N+2 s}} d x d y.
$$

We now state our basic assumptions on $V(x)$, $f(x)$ and $g(x)$.
\begin{itemize}
	\item[$(V)$] Assume $V(x)\in L^\infty(\mathbb{R}^N)$, and 
	$$
	0\leq   V(x) \leq V_{\infty}=\lim _{|x| \rightarrow \infty} V(x), \quad \text { for all } x \in \mathbb{R}^N.
	$$
	\item[$(f)$] Assume $f(x)\in L^\infty(\mathbb{R}^N)$, and 
	$$
	0< f_{\infty}=\lim _{|x| \rightarrow \infty} f(x) \leq f(x) \leq f_{max}, \quad \text { for all } x \in \mathbb{R}^N.
	$$
	\item[$(g)$] Assume $g(x)\in L^\infty(\mathbb{R}^N)$ , and 
	$$
	0< g_{\infty}=\lim _{|x| \rightarrow \infty} g(x) \leq g(x) \leq g_{max}, \quad \text { for all } x \in \mathbb{R}^N.
	$$
\end{itemize}

To study the normalized solutions of  \eqref{eq1}, we analyze the critical points of the energy functional:
\begin{equation}\label{energyfunctional}
	\begin{aligned}
		J(u)=&\frac{1}{2}\int_{\mathbb{R}^N}|(-\Delta)^{\frac{s}{2}} u|^2dx+\frac{1}{2}\int_{\mathbb{R}^N}V(x)u^2dx-\frac{1}{2q} \int_{\mathbb{R}^N}\left(I_\alpha *(f|u|^q)\right)f|u|^q d x\\
		&-\frac{1}{2p}\int_{\mathbb{R}^N}\left(I_\alpha *(g|u|^p)\right)g|u|^pd x,
	\end{aligned}
\end{equation}
subject to the constraint:
$$
S(a)=\left\{u \in H^s\left(\mathbb{R}^N\right): \int_{\mathbb{R}^N}|u|^2 d x=a\right\}.
$$

\begin{lemma}\cite{lieb-loss} \label{le1.1} Let $r, t>1$ and $\alpha \in(0, N)$ with $\frac{1}{r}+\frac{1}{t}=1+\frac{\alpha}{N}$. For $f \in L^r\left(\mathbb{R}^N\right)$ and $h \in L^t\left(\mathbb{R}^N\right)$,  there exists a sharp constant $C(r, t, \alpha, N)$ independent of $f$ and $h$, such that
	$$
	\int_{\mathbb{R}^N} \int_{\mathbb{R}^N} \frac{f(x) h(y)}{|x-y|^{N-\alpha}} d x d y \leq C(r, t, \alpha, N)\|f\|_r\|h\|_t .
	$$
\end{lemma}

By this lemma, the functional $J$ defined in \eqref{energyfunctional} is well-defined on  $H^s\left(\mathbb{R}^N\right)$ and belongs to $ C^1$. For convenience, we define
$$Q(u)=\int_{\mathbb{R}^N}\left(I_\alpha *(f|u|^q)\right)f|u|^q d x$$  and $$P(u)=\int_{\mathbb{R}^N}\left(I_\alpha *(g|u|^p)\right)g|u|^pd x.$$ 
Thus,  $J$ can be rewritten as:
\[J(u)=\frac{1}{2}\|u\|^2+\frac{1}{2}\int_{\mathbb{R}^N}V(x)u^2dx-\frac{1}{2q} Q(u)-\frac{1}{2p}P(u).\]

\begin{lemma}\cite{feng-zhang}\label{le1.2}  Let $N>2 s, 0<s<1$ and $\frac{N+\alpha}{N}<t<2_{\alpha, s}^*$. Then, for all $u \in H^s\left(\mathbb{R}^N\right)$,
	$$
	\int_{\mathbb{R}^N}\left(I_\alpha *|u|^t\right)|u|^t d x \leq C_{\alpha, t}\left\|u\right\|^{2 t \gamma_{t, s}}\|u\|_2^{2 t\left(1-\gamma_{t, s}\right)},
	$$
	where $\gamma_{t, s}=\frac{N t-N-\alpha}{2 t s}$ and the optimal constant $C_{\alpha, t}$ is given by
	$$
	C_{\alpha, t}=\frac{2 s t}{2 s t-N t+N+\alpha}\left(\frac{2 s t-N t+N+\alpha}{N t-N-\alpha}\right)^{\frac{Nt-N-\alpha}{2 s}}\|U\|_2^{2-2 t} \text {, }
	$$
	where $U$ is the ground state solution of
	$$
	(-\Delta)^s U+U-\left(I_\alpha *|U|^t\right)|U|^{t-2} U=0 .
	$$
\end{lemma}

\begin{lemma}\cite{lieb,Moroz}\label{le1.3}
	For all $u \in H^s\left(\mathbb{R}^N\right)$, 
	$$
	\int_{\mathbb{R}^N}\left(I_\alpha *|u|^{1+\frac{\alpha}{N}}\right)|u|^{1+\frac{\alpha}{N}} d x \leq S_\alpha^{-\left(1+\frac{\alpha}{N}\right)}\left(\int_{\mathbb{R}^N}|u|^2 d x\right)^{1+\frac{\alpha}{N}},
	$$
	where $S_\alpha$ is the infimum of the minimization problem
	$$
	S_\alpha=\inf \left\{\int_{\mathbb{R}^N}|u|^2 d x: u \in L^2\left(\mathbb{R}^N\right) \text { and } \int_{\mathbb{R}^N}\left(I_\alpha *|u|^{1+\frac{\alpha}{N}}\right)|u|^{1+\frac{\alpha}{N}} d x=1\right\}>0.
	$$
	Moreover,  $S_\alpha$ is achieved if and only if
	\[ u(x)=V_\epsilon(x)=K\left(\frac{\epsilon}{\epsilon^2+|x-z|^2}\right)^{\frac{N}{2}}\]
	for some  $K \in \mathbb{R}$, $\epsilon \in \mathbb{R}^{+}$ and $z \in \mathbb{R}^N$.
\end{lemma}

Set $\gamma_{\frac{N+\alpha}{N}, s}=0$ and let $C_{\alpha,\frac{N+\alpha}{N}}=S_\alpha^{-\left(1+\frac{\alpha}{N}\right)}$. Combining Lemmas \ref{le1.2} and \ref{le1.3}, we  derive the  inequality:
\begin{equation}\label{1.3}
	\int_{\mathbb{R}^N}\left(I_\alpha *|u|^t\right)|u|^t d x \leq C_{\alpha, t}\left\| u\right\|^{2 t \gamma_{t, s}}\|u\|_2^{2 t\left(1-\gamma_{t, s}\right)},
\end{equation}
where $t\in\left[\frac{N+\alpha}{N},2_{\alpha,s}^*\right)$.

\vspace{10pt}

From the definition of $\gamma_{t,s}$ in Lemma \ref{le1.2},  by direct computation, we have 
$$
t \gamma_{t,s}\left\{\begin{array}{ll}
	=0, & t=\frac{N+\alpha}{N}, \\
	<1, & \frac{N+\alpha}{N}<t<\frac{N+\alpha+2s}{N}, \\
	=1, & t=\frac{N+\alpha+2s}{N}.
\end{array}\right.
$$

We now can present our main results.

\begin{theorem}\label{thm1}
	Assume $\frac{N+\alpha}{N}= q<p< \frac{N+2 s+\alpha}{N}$, and that $(V)$, $(f)$ and $(g)$ hold. Then for any $a>0$, there exists $V_*>0$ such that for  $\|V\|_{\infty}<V_*$,  equation \eqref{eq1} admits at least a  weak solution $\left(u, \lambda\right) \in H^s\left(\mathbb{R}^N\right) \times \mathbb{R}$  satisfying $\|u\|_2^2=a$, $\lambda<-f^2_{\infty}S_\alpha^{-\frac{N+\alpha}{N}} a^{\frac{\alpha}{N}}$, and $J(u)<-\frac{f^2_{\infty}N}{2(N+\alpha)} S_\alpha^{-\frac{N+\alpha}{N}} a^{\frac{N+\alpha}{N}}$.
\end{theorem}

\begin{theorem}\label{thm2}
	Assume $\frac{N+\alpha}{N}<q<p=\frac{N+2 s+\alpha}{N}$, and that $(V)$, $(f)$ and $(g)$ hold. Then for any $0<a<a^*=\left(\frac{p}{C_{\alpha, p}g^2_{max}}\right)^{\frac{1}{p-1}}$, there exists $V_*>0$ such that for  $\|V\|_{\infty}<V_*$,  equation \eqref{eq1} admits at least a  weak solution $\left(u, \lambda\right) \in H^s\left(\mathbb{R}^N\right) \times \mathbb{R}$ satisfying $\|u\|_2^2=a$, $\lambda<0$, and $J(u)<0$.
\end{theorem}

In this paper, we can also consider  the existence of  normalized solutions for the following nonlinear fractional Choquard equation:
\begin{equation}\label{neq1}
	\begin{aligned}
		(-\Delta)^s u+V(x)u=&\lambda u+f(x)\left(I_\alpha *(f|u|^q)\right)|u|^{q-2} u\\
		&+\left(I_\alpha *|u|^k\right)|u|^{k-2} u+g(x)\left(I_\alpha *(g|u|^p)\right)|u|^{p-2} u, \quad x \in \mathbb{R}^N
	\end{aligned}
\end{equation}
subject to the constraint
$$
\int_{\mathbb{R}^N}|u|^2 d x=a>0.
$$
Here $N>2 s, s \in(0,1), \alpha \in(0, N)$,  $\frac{N+\alpha}{N}=q<k<p=\frac{N+2 s+\alpha}{N}$, and $\lambda \in \mathbb{R}$  is an unknown parameter acting as a Lagrange multiplier.

\begin{theorem}\label{thm3}
	Assume $\frac{N+\alpha}{N}=q<k<p=\frac{N+2 s+\alpha}{N}$, and that $(V)$, $(f)$ and $(g)$ hold. Then for any $0<a<a^*=\left(\frac{p}{C_{\alpha, p}g^2_{max}}\right)^{\frac{1}{p-1}}$, there exists $V_*>0$ such that for  $\|V\|_{\infty}<V_*$, equation \eqref{neq1} admits at least a  weak solution $\left(u, \lambda\right) \in H^s\left(\mathbb{R}^N\right) \times \mathbb{R}$   satisfying $\|u\|_2^2=a$, $\lambda<-f^2_{\infty}S_\alpha^{-\frac{N+\alpha}{N}} a^{\frac{\alpha}{N}}$, and $I(u)<-\frac{f^2_{\infty}N}{2(N+\alpha)}
	S_\alpha^{-\frac{N+\alpha}{N}} a^{\frac{N+\alpha}{N}}$.
\end{theorem}

The remainder of this paper is organized as follows. In Section 2, we establish the variational preliminaries needed for the minimization procedure on the mass constraint. In Section 3, we complete the proofs of Theorems \ref{thm1}, \ref{thm2}, and \ref{thm3} by applying the compactness result obtained in Section 2.

\section{Preliminaries}

\begin{lemma}\label{le2.1}
	\begin{itemize}
		\item[$(1)$]  If $\frac{N+\alpha}{N}=q<p<\frac{N+\alpha+2s}{N}$,  then for any $a>0$, the functional $J$ is bounded from below on the set $S\left(a\right)$.
		\item[$(2)$] If  $\frac{N+\alpha}{N}<q<p=\frac{N+\alpha+2s}{N}$, then for any $0<a<a^*=\left(\frac{p}{C_{\alpha, p}g^2_{max}}\right)^{\frac{1}{p-1}}$, the functional $J$ is bounded from below on the set  $S\left(a\right)$.
	\end{itemize}
\end{lemma} 
\begin{proof}
	Consider the constrained  functional $J|_{S(a)}$. For every $u \in S(a)$, by \eqref{1.3}, we have
	\begin{equation}\label{2.1}
		\begin{aligned}
			J(u) &\geq \frac{1}{2}\|u\|^2-\frac{f^2_{max}}{2 q} C_{\alpha, q}\|u\|^{2q \gamma_{q, s}}\|u\|_2^{2 q\left(1-\gamma_{q, s}\right)}-\frac{g^2_{max}}{2 p} C_{\alpha, p}\|u\|^{2 p \gamma_{p, s}}\|u\|_2^{2 p\left(1-\gamma_{p, s}\right)} \\
			& =\frac{1}{2}\|u\|^2-\frac{ f^2_{max}a^{q\left(1-\gamma_{q,s}\right)}}{2q} C_{\alpha, q}\|u\|^{2 q \gamma_{q, s}}-\frac{g^2_{max}a^{p\left(1-\gamma_{p, s}\right)}}{2p} C_{\alpha, p}\|u\|^{2 p \gamma_{p, s}}.
		\end{aligned}
	\end{equation}
	
	Case 1. $\frac{N+\alpha}{N}=q<p<\frac{N+\alpha+2s}{N}$.  
	Since $0=2q \gamma_{q, s}<2 p \gamma_{p, s}<2$, it follows from \eqref{2.1} that $J$ is bounded from below and coercive on $S(a)$.

	Case 2. $\frac{N+\alpha}{N}<q<p=\frac{N+\alpha+2s}{N}$. Since $p \gamma_{p, s}=1$, it follows from  \eqref{2.1} that
	\begin{equation*}
		J(u) \geq \left[\frac{1}{2}-\frac{g^2_{max}a^{p-1}}{2p} C_{\alpha, p}\right]\|u\|^2-\frac{ f^2_{max}a^{q\left(1-\gamma_{q,s}\right)}}{2q} C_{\alpha, q}\|u\|^{2 q \gamma_{q, s}}.
	\end{equation*}
	Since $a<a^*$, the term $\frac{1}{2}-\frac{g^2_{max}a^{p-1}}{2p} C_{\alpha, p}$ is positive. Moreover, as $2q \gamma_{q, s}<2$, we conclude that $J$ is bounded from below and coercive on $S(a)$.  
\end{proof}

\begin{lemma}\label{le2.2} 
	Let \(a>0\) in case \((1)\), or \(0<a<a^*\) in case \((2)\). Then there exists a constant \(V_*=V_*(a)>0\) such that, for all \(\|V\|_{\infty}\in[0,V_*)\), the following hold:
	\begin{itemize}
		\item[$(1)$] If \(\frac{N+\alpha}{N}=q<p<\frac{N+\alpha+2s}{N}\), then
		\[
		m(a)<-\frac{f^2_{\infty}N}{2(N+\alpha)}
		S_\alpha^{-\frac{N+\alpha}{N}}a^{\frac{N+\alpha}{N}}.
		\]
		\item[$(2)$] If \(\frac{N+\alpha}{N}<q<p=\frac{N+\alpha+2s}{N}\), then
		\[
		m(a)<0.
		\]
	\end{itemize}
\end{lemma}
\begin{proof}
	By Lemma \ref{le2.1}, $m(a)$ is well-defined. From Lemma \ref{le1.3},  there exists $v$  satisfying  
	\begin{equation}\label{n2.2}
		S_\alpha\left(\int_{\mathbb{R}^N}\left(I_\alpha *|v|^{\frac{N+\alpha}{N}}\right)|v|^{\frac{N+\alpha}{N}} d x\right)^{\frac{N}{N+\alpha}}=\int_{\mathbb{R}^N}|v|^2 d x.
	\end{equation}
	Define $u=\frac{\sqrt{a} v}{\|v\|_2}$ and  consider the rescaled function $u_t(x)=t^{\frac{N}{2}} u(t x)$. Note that $u_t \in S\left(a\right)$ for all $t>0$. A direct computation yields
	\begin{equation}\label{n2.3}
		J\left(u_t\right) \leq \frac{1}{2} t^{2s} \|u\|^2+\frac{\|V\|_\infty}{2}a-\frac{f^2_{\infty}}{2q} t^{Nq-N-\alpha} \widetilde{Q}(u)-\frac{g^2_{\infty}}{2p} t^{Np-N-\alpha} \widetilde{P}(u),
	\end{equation}
	where \[\widetilde{Q}(u)=\int_{\mathbb{R}^N}\left(I_\alpha *|u|^q\right)|u|^q d x \] and \[\widetilde{P}(u)=\int_{\mathbb{R}^N}\left(I_\alpha *|u|^p\right)|u|^p d x .\]

	Case 1. $\frac{N+\alpha}{N}=q<p<\frac{N+\alpha+2s}{N}$.  From \eqref{n2.2},  we have $\widetilde{Q}(u)=a^qS_\alpha^{-q}$. Substituting this into \eqref{n2.3} gives
	\begin{equation}\label{n2.4}
		J\left(u_t\right) \leq \frac{1}{2} t^{2s} \|u\|^2+\frac{\|V\|_\infty}{2}a-\frac{f^2_{\infty}}{2q} a^qS_\alpha^{-q}-\frac{g^2_{\infty}}{2p} t^{Np-N-\alpha} \widetilde{P}(u).
	\end{equation}
	Since \(p < \frac{N + 2s + \alpha}{N}\), the dominant term for small $t > 0$ is \(t^{Np-N-\alpha}\), and thus   \[\frac{1}{2} t^{2s} \|u\|^2-\frac{g^2_{\infty}}{2p} t^{Np-N-\alpha} \widetilde{P}(u)<0\] for \(t\) sufficiently small. Consequently, there exists $V_*>0$ such that for $\|V\|_\infty \in [0, V_*)$,
	\[\frac{1}{2} t^{2s} \|u\|^2-\frac{g^2_{\infty}}{2p} t^{Np-N-\alpha} \widetilde{P}(u)+\frac{\|V\|_\infty}{2}a<0.\]
	Combining this with \eqref{n2.4}, we obtain
	\begin{equation*}
		J\left(u_t\right) <-\frac{f^2_{\infty}}{2q} a^qS_\alpha^{-q}=-\frac{f^2_{\infty}N}{2(N+\alpha)} S_\alpha^{-\frac{N+\alpha}{N}} a^{\frac{N+\alpha}{N}}.
	\end{equation*}

	Case 2. $\frac{N+\alpha}{N}<q<p=\frac{N+\alpha+2s}{N}$. Here, $q < \frac{N + 2s + \alpha}{N}$, so for small \(t > 0\),    \[\frac{1}{2} t^{2s} \|u\|^2-\frac{f^2_{\infty}}{2q} t^{Nq-N-\alpha} \widetilde{Q}(u)<0.\] Thus,  there exists  $V_*>0$ such that for $\|V\|_\infty \in [0, V_*)$,
	\[\frac{1}{2} t^{2s} \|u\|^2-\frac{f^2_{\infty}}{2q} t^{Nq-N-\alpha} \widetilde{Q}(u)+\frac{\|V\|_\infty}{2}a<0.\]
	From \eqref{n2.3}, this implies
	\begin{equation*}
		J\left(u_t\right) <0,
	\end{equation*}
	completing the proof.  
\end{proof}

\begin{lemma}\label{le2.3}  
	The following hold:
	\begin{itemize}
		\item[$(1)$] If $\frac{N+\alpha}{N}=q<p<\frac{N+\alpha+2s}{N}$, then the function \(m(a)\) is continuous on \((0,+\infty)\).
		\item[$(2)$] If $\frac{N+\alpha}{N}<q<p=\frac{N+\alpha+2s}{N}$, then the function \(m(a)\) is continuous on \((0,a^*)\).
	\end{itemize}
\end{lemma}
\begin{proof}
	It is enough to prove sequential continuity. Let \(b\in(0,+\infty)\) in case \((1)\), or \(b\in(0,a^*)\) in case \((2)\). Let \(\{a_n\}\) be a sequence in the corresponding interval such that
	\[
	a_n\to b \qquad \text{as } n\to\infty .
	\]
	For each \(n\in\mathbb N\), choose \(u_n\in S(a_n)\) such that
	\begin{equation}\label{2.6}
		J(u_n)<m(a_n)+\frac{1}{n}.
	\end{equation}
	
	We first note that \(\{u_n\}\) is bounded in \(H^s(\mathbb R^N)\). Indeed, since \(a_n\to b\), the sequence \(\{a_n\}\) remains in a compact subinterval of the corresponding admissible interval. Taking a fixed \(z\in S(b)\) and setting
	\[
	z_n=\sqrt{\frac{a_n}{b}}z\in S(a_n),
	\]
	we have
	\[
	m(a_n)\leq J(z_n)\leq C .
	\]
	Combining this estimate with \eqref{2.6} and the lower bounds obtained in Lemma \ref{le2.1}, we obtain the boundedness of \(\{u_n\}\) in \(H^s(\mathbb R^N)\). Hence, by \eqref{1.3}, the quantities \(Q(u_n)\) and \(P(u_n)\) are also bounded.
	
	Define
	\[
	v_n=\sqrt{\frac{b}{a_n}}u_n .
	\]
	Then \(v_n\in S(b)\). Since
	\[
	Q(v_n)=\left(\frac{b}{a_n}\right)^q Q(u_n),
	\qquad
	P(v_n)=\left(\frac{b}{a_n}\right)^p P(u_n),
	\]
	we have
	\[
	\begin{aligned}
		J(v_n)
		=&J(u_n)
		+\frac12\left(\frac{b}{a_n}-1\right)\|u_n\|^2
		+\frac12\left(\frac{b}{a_n}-1\right)
		\int_{\mathbb R^N}V(x)|u_n|^2\,dx \\
		&-\frac{1}{2q}
		\left[\left(\frac{b}{a_n}\right)^q-1\right]Q(u_n)
		-\frac{1}{2p}
		\left[\left(\frac{b}{a_n}\right)^p-1\right]P(u_n).
	\end{aligned}
	\]
	Since \(a_n\to b\), and since \(\{u_n\}\), \(\{Q(u_n)\}\), and \(\{P(u_n)\}\) are bounded, it follows that
	\[
	J(v_n)=J(u_n)+o_n(1).
	\]
	Therefore, by \eqref{2.6},
	\[
	m(b)\leq J(v_n)\leq m(a_n)+\frac{1}{n}+o_n(1).
	\]
	Consequently,
	\begin{equation}\label{2.9}
		m(b)\leq \liminf_{n\to\infty}m(a_n).
	\end{equation}
	
	Conversely, let \(\{e_n\}\subset S(b)\) be a minimizing sequence for \(m(b)\), namely
	\[
	J(e_n)=m(b)+o_n(1).
	\]
	By Lemma \ref{le2.1}, the sequence \(\{e_n\}\) is bounded in \(H^s(\mathbb R^N)\). Define
	\[
	f_n=\sqrt{\frac{a_n}{b}}e_n .
	\]
	Then \(f_n\in S(a_n)\). Using the same scaling computation as above, and using \(a_n\to b\), we obtain
	\[
	J(f_n)=J(e_n)+o_n(1).
	\]
	Thus
	\[
	m(a_n)\leq J(f_n)=J(e_n)+o_n(1)=m(b)+o_n(1),
	\]
	which gives
	\begin{equation}\label{2.10}
		\limsup_{n\to\infty}m(a_n)\leq m(b).
	\end{equation}
	Combining \eqref{2.9} with \eqref{2.10}, we conclude that
	\[
	\lim_{n\to\infty}m(a_n)=m(b).
	\]
	Hence \(m\) is continuous on the corresponding interval. The proof is complete.
\end{proof}

\begin{lemma}\label{le2.4}  
	Let \(a_2\) be admissible, namely \(a_2>0\) if
	\[
	\frac{N+\alpha}{N}=q<p<\frac{N+\alpha+2s}{N},
	\]
	and \(0<a_2<a^*\) if
	\[
	\frac{N+\alpha}{N}<q<p=\frac{N+\alpha+2s}{N}.
	\]
	Assume that \(\|V\|_\infty\in[0,V_*)\), where \(V_*=V_*(a_2)>0\) is chosen so that
	\[
	m(a_2)<0.
	\]
	Then, for every admissible \(a_1\) with \(0<a_1<a_2\), one has
	\[
	\frac{a_1}{a_2}m(a_2)<m(a_1).
	\]
\end{lemma} 
\begin{proof}
	Let
	\[
	\xi=\frac{a_2}{a_1}>1.
	\]
	If \(m(a_1)\geq0\), then, since \(m(a_2)<0\), we immediately have
	\[
	\frac{a_1}{a_2}m(a_2)<0\leq m(a_1).
	\]
	Thus it remains to consider the case
	\[
	m(a_1)<0.
	\]
	
	Let \(\{u_n\}\subset S(a_1)\) be a minimizing sequence for \(m(a_1)\), namely
	\[
	J(u_n)\to m(a_1)\quad\text{as }n\to+\infty .
	\]
	Define
	\[
	v_n=\sqrt{\xi}\,u_n .
	\]
	Then \(v_n\in S(a_2)\). Moreover,
	\[
	Q(v_n)=\xi^q Q(u_n),
	\qquad
	P(v_n)=\xi^p P(u_n).
	\]
	Hence
	\begin{equation}\label{n2.8}
		\begin{aligned}
			m(a_2)
			\leq J(v_n)
			=&\xi J(u_n)
			-\frac{\xi^q-\xi}{2q}Q(u_n)
			-\frac{\xi^p-\xi}{2p}P(u_n).
		\end{aligned}
	\end{equation}
	Since \(\xi>1\) and \(q,p>1\), we have
	\[
	\frac{\xi^q-\xi}{2q}>0,
	\qquad
	\frac{\xi^p-\xi}{2p}>0.
	\]
	
	We claim that there exist \(C>0\) and \(n_0\in\mathbb N\) such that
	\[
	Q(u_n)+P(u_n)\geq C,
	\qquad n\geq n_0 .
	\]
	Indeed, otherwise there exists a subsequence, still denoted by \(\{u_n\}\), such that
	\[
	Q(u_n)+P(u_n)\to0.
	\]
	Since \(Q(u_n)\geq0\) and \(P(u_n)\geq0\), this gives
	\[
	Q(u_n)\to0,
	\qquad
	P(u_n)\to0.
	\]
	Using \(V(x)\geq0\), we obtain
	\[
	J(u_n)
	=
	\frac12\|u_n\|^2
	+\frac12\int_{\mathbb R^N}V(x)|u_n|^2\,dx
	-\frac{1}{2q}Q(u_n)
	-\frac{1}{2p}P(u_n)
	\geq o_n(1).
	\]
	Letting \(n\to+\infty\), we get
	\[
	m(a_1)\geq0,
	\]
	which contradicts \(m(a_1)<0\). Therefore the claim holds.
	
	Set
	\[
	\eta=\min\left\{
	\frac{\xi^q-\xi}{2q},
	\frac{\xi^p-\xi}{2p}
	\right\}>0.
	\]
	From \eqref{n2.8}, for all \(n\geq n_0\),
	\[
	m(a_2)
	\leq \xi J(u_n)-\eta\bigl(Q(u_n)+P(u_n)\bigr)
	\leq \xi J(u_n)-\eta C .
	\]
	Passing to the limit as \(n\to+\infty\), we obtain
	\[
	m(a_2)\leq \xi m(a_1)-\eta C<\xi m(a_1).
	\]
	Since \(\xi=\frac{a_2}{a_1}\), this is equivalent to
	\[
	\frac{a_1}{a_2}m(a_2)<m(a_1).
	\]
	The proof is complete.
\end{proof}

Since we consider equation \eqref{eq1} in the whole space $\mathbb{R}^N$, the loss of compactness must be addressed. To overcome this, we introduce the following energy functional:
\begin{equation}\label{energyfunctional2}
	\begin{aligned}
		J_\infty(u)=&\frac{1}{2}\int_{\mathbb{R}^N}|(-\Delta)^{\frac{s}{2}} u|^2dx+\frac{V_\infty}{2}\int_{\mathbb{R}^N} u^2dx-\frac{f^2_\infty}{2q} \int_{\mathbb{R}^N}\left(I_\alpha *|u|^q\right)|u|^q d x\\
		&-\frac{g^2_\infty}{2p}\int_{\mathbb{R}^N}\left(I_\alpha *|u|^p\right)|u|^pd x,
	\end{aligned}
\end{equation}
and define
\[m_\infty(a)=\inf\limits _{u \in S\left(a\right)} J_\infty(u).\]
Using similar arguments, one can show that  $J_\infty$ and $m_\infty(a)$ satisfy the same properties as those stated in  Lemma \ref{le2.1}
through Lemma \ref{le2.4}.

\begin{lemma}\label{le2.5}
	\begin{itemize}
		\item[$(1)$] If $\frac{N+\alpha}{N}=q<p<\frac{N+\alpha+2s}{N}$, then for any $a>0$, $m(a) \leq m_\infty(a)$.
		\item[$(2)$] If $\frac{N+\alpha}{N}<q<p=\frac{N+\alpha+2s}{N}$, then for any $0<a <a^*$, we have $m(a) \leq m_\infty(a)$.
	\end{itemize}
\end{lemma}
\begin{proof}
	For \(a > 0\) (or \(0 < a < a^*\) in case $(2)$), let \(\{u_n\} \subset S(a)\) be a minimizing sequence for \(m_\infty(a)\), i.e.,  
	\[
	J_\infty(u_n) \to m_\infty(a) \quad \text{as } n \to +\infty.
	\]  
	Since \(V(x) \leq V_\infty\), \(f(x) \geq f_\infty\), and \(g(x) \geq g_\infty\), we have  
	\[
	m(a) \leq J(u_n) \leq J_\infty(u_n) = m_\infty(a) + o_n(1).
	\]  
	Taking \(n \to +\infty\), the result follows. 
\end{proof}

\begin{lemma}\label{le2.6}
	Let \(a\) be admissible, namely \(a\in(0,+\infty)\) if
	\[
	\frac{N+\alpha}{N}=q<p<\frac{N+\alpha+2s}{N},
	\]
	and \(a\in(0,a^*)\) if
	\[
	\frac{N+\alpha}{N}<q<p=\frac{N+\alpha+2s}{N}.
	\]
	Let \(V_*=V_*(a)>0\) be chosen as in Lemma \ref{le2.2}, and assume that
	\[
	\|V\|_\infty\in[0,V_*).
	\]
	If \(\{u_n\}\subset S(a)\) is a minimizing sequence for \(m(a)\), then, up to a subsequence, one of the following holds:
	\begin{itemize}
		\item[$(I)$] $\left\{u_n\right\}$  converges strongly in $H^s\left(\mathbb{R}^N\right)$.
		\item[$(II)$] There exists $\left\{y_n\right\} \subset \mathbb{R}^N$ with $\left|y_n\right| \rightarrow \infty$ such that  $v_n(x)=u_n\left(x+y_n\right)$  converges strongly to a function $v \in S\left(a\right)$ satisfying $J(v)=m(a)$.
	\end{itemize}
\end{lemma} 
\begin{proof}
	By Lemmas \ref{le2.1} and \ref{le2.2}, \(\{u_n\}\) is bounded in
	\(H^s(\mathbb R^N)\). Hence, up to a subsequence,
	\[
	u_n\rightharpoonup u \quad \text{in } H^s(\mathbb R^N),
	\qquad
	u_n\to u \quad \text{a.e. in } \mathbb R^N .
	\]
	We divide the proof into two cases.
	
	{\bf Case \(1\): \(u\not\equiv0\).} Assume first that
	\[
	\|u\|_2^2=b\neq a .
	\]
	By the weak lower semicontinuity of the \(L^2\)-norm, \(b\in(0,a)\). Set
	\[
	w_n=u_n-u .
	\]
	Then, by the Brezis--Lieb lemma and the nonlocal Brezis--Lieb splitting, see Lemma 3.6 of \cite{meng-he},
	\begin{equation}\label{2.12}
		\|u_n\|^2=\|w_n\|^2+\|u\|^2+o_n(1),
	\end{equation}
	\begin{equation}\label{2.13}
		\int_{\mathbb R^N}V(x)|u_n|^2\,dx
		=
		\int_{\mathbb R^N}V(x)|w_n|^2\,dx
		+
		\int_{\mathbb R^N}V(x)|u|^2\,dx
		+o_n(1),
	\end{equation}
	and
	\begin{equation}\label{2.14}
		Q(u_n)=Q(w_n)+Q(u)+o_n(1),
		\qquad
		P(u_n)=P(w_n)+P(u)+o_n(1).
	\end{equation}
	Let
	\[
	d_n=\|w_n\|_2^2 .
	\]
	Then \(d_n\to d=a-b\). In particular, \(d\in(0,a)\), and \(d_n\in(0,a)\) for \(n\) large enough. From \eqref{2.12}--\eqref{2.14}, we obtain
	\[
	\begin{aligned}
		m(a)+o_n(1)
		=J(u_n)
		=&\frac12\|w_n\|^2
		+\frac12\int_{\mathbb R^N}V(x)|w_n|^2\,dx
		-\frac{1}{2q}Q(w_n)
		-\frac{1}{2p}P(w_n)  \\
		&+\frac12\|u\|^2
		+\frac12\int_{\mathbb R^N}V(x)|u|^2\,dx
		-\frac{1}{2q}Q(u)
		-\frac{1}{2p}P(u)
		+o_n(1)  \\
		=&J(w_n)+J(u)+o_n(1)  \\
		\geq& m(d_n)+m(b)+o_n(1).
	\end{aligned}
	\]
	Letting \(n\to\infty\) and using Lemma \ref{le2.3}, we get
	\begin{equation}\label{2.16}
		m(a)\geq m(d)+m(b).
	\end{equation}
	By Lemma \ref{le2.4}, since \(b,d\in(0,a)\) and \(b+d=a\), we have
	\[
	m(d)+m(b)>
	\frac{d}{a}m(a)+\frac{b}{a}m(a)=m(a),
	\]
	which contradicts \eqref{2.16}. Therefore,
	\[
	\|u\|_2^2=a .
	\]
	Since \(u_n\rightharpoonup u\) in \(L^2(\mathbb R^N)\) and \(\|u_n\|_2\to\|u\|_2\), we have
	\[
	u_n\to u \quad \text{in } L^2(\mathbb R^N).
	\]
	By interpolation,
	\begin{equation}\label{nn2.16}
		u_n\to u \quad \text{in } L^t(\mathbb R^N),
		\qquad 2<t<2_s^* .
	\end{equation}
	Together with \eqref{1.3}, this gives
	\[
	Q(u_n)\to Q(u),
	\qquad
	P(u_n)\to P(u).
	\]
	Using the weak lower semicontinuity of the norm and of the potential term, we deduce that
	\[
	\begin{aligned}
		m(a)
		=&\lim_{n\to\infty}J(u_n)  \\
		=&\lim_{n\to\infty}
		\left[
		\frac12\|u_n\|^2
		+\frac12\int_{\mathbb R^N}V(x)|u_n|^2\,dx
		-\frac{1}{2q}Q(u_n)
		-\frac{1}{2p}P(u_n)
		\right]  \\
		\geq& J(u)\geq m(a).
	\end{aligned}
	\]
	Hence \(J(u)=m(a)\), and
	\[
	\|u_n\|\to\|u\|.
	\]
	Since \(H^s(\mathbb R^N)\) is a Hilbert space, we conclude that
	\[
	u_n\to u \quad \text{in } H^s(\mathbb R^N).
	\]
	Thus alternative \((I)\) holds.
	
	{\bf Case \(2\): \(u\equiv0\).} Then
	\[
	u_n\rightharpoonup0 \quad \text{in } H^s(\mathbb R^N),
	\qquad
	u_n\to0 \quad \text{in } L^2_{\mathrm{loc}}(\mathbb R^N).
	\]
	Since
	\[
	\lim_{|x|\to\infty}V(x)=V_\infty ,
	\]
	for every \(\delta>0\), there exists \(R=R(\delta)>0\) such that
	\[
	|V(x)-V_\infty|<\delta
	\qquad \text{for } |x|>R .
	\]
	Therefore,
	\[
	\begin{aligned}
		\left|
		\int_{\mathbb R^N}(V(x)-V_\infty)|u_n|^2\,dx
		\right|
		&\leq C\int_{B_R(0)}|u_n|^2\,dx
		+\delta\int_{B_R^c(0)}|u_n|^2\,dx  \\
		&\leq o_n(1)+\delta a .
	\end{aligned}
	\]
	Letting \(\delta\to0\), we get
	\begin{equation}\label{n2.15}
		\int_{\mathbb R^N}V(x)|u_n|^2\,dx
		=
		V_\infty\int_{\mathbb R^N}|u_n|^2\,dx
		+o_n(1).
	\end{equation}
	
	We next compare the nonlocal terms with their limiting counterparts. Since
	\[
	\lim_{|x|\to\infty}f(x)=f_\infty,
	\qquad
	\lim_{|x|\to\infty}g(x)=g_\infty,
	\]
	and since \(u_n\to0\) locally, the Hardy--Littlewood--Sobolev inequality together with \eqref{1.3} yields
	\begin{equation}\label{n2.16}
		\int_{\mathbb R^N}
		\left(I_\alpha *[(f-f_\infty)|u_n|^q]\right)
		(f-f_\infty)|u_n|^q\,dx=o_n(1),
	\end{equation}
	\begin{equation}\label{n2.17}
		\int_{\mathbb R^N}
		\left(I_\alpha *[(f-f_\infty)|u_n|^q]\right)
		f_\infty |u_n|^q\,dx=o_n(1),
	\end{equation}
	and
	\begin{equation}\label{n2.18}
		\int_{\mathbb R^N}
		\left(I_\alpha *(f_\infty |u_n|^q)\right)
		(f-f_\infty)|u_n|^q\,dx=o_n(1).
	\end{equation}
	Consequently,
	\begin{equation}\label{n2.19}
		Q(u_n)
		=
		f_\infty^2
		\int_{\mathbb R^N}
		\left(I_\alpha*|u_n|^q\right)|u_n|^q\,dx
		+o_n(1).
	\end{equation}
	Similarly,
	\begin{equation}\label{n2.20}
		P(u_n)
		=
		g_\infty^2
		\int_{\mathbb R^N}
		\left(I_\alpha*|u_n|^p\right)|u_n|^p\,dx
		+o_n(1).
	\end{equation}
	From \eqref{n2.15}, \eqref{n2.19}, and \eqref{n2.20}, it follows that
	\[
	J_\infty(u_n)=J(u_n)+o_n(1).
	\]
	Hence
	\[
	J_\infty(u_n)\to m(a).
	\]
	
	We claim that there exist \(R>0\), \(\beta>0\), and \(\{y_n\}\subset\mathbb R^N\) such that
	\begin{equation}\label{nne2.18}
		\int_{B_R(y_n)}|u_n|^2\,dx\geq\beta
		\qquad \text{for all } n .
	\end{equation}
	Indeed, if this were false, then by the fractional Lions vanishing lemma,
	\[
	u_n\to0 \quad \text{in } L^t(\mathbb R^N),
	\qquad 2<t<2_s^* .
	\]
	
	If
	\[
	\frac{N+\alpha}{N}=q<p<\frac{N+\alpha+2s}{N},
	\]
	then the \(p\)-term vanishes, and the lower Hardy--Littlewood--Sobolev inequality gives
	\[
	\begin{aligned}
		m(a)+o_n(1)
		=J_\infty(u_n)
		&\geq
		\frac12\|u_n\|^2
		-\frac{f_\infty^2}{2q}
		\int_{\mathbb R^N}
		\left(I_\alpha*|u_n|^q\right)|u_n|^q\,dx
		+o_n(1)  \\
		&\geq
		\frac12\|u_n\|^2
		-\frac{f_\infty^2N}{2(N+\alpha)}
		S_\alpha^{-\frac{N+\alpha}{N}}
		a^{\frac{N+\alpha}{N}}
		+o_n(1)  \\
		&\geq
		-\frac{f_\infty^2N}{2(N+\alpha)}
		S_\alpha^{-\frac{N+\alpha}{N}}
		a^{\frac{N+\alpha}{N}}
		+o_n(1),
	\end{aligned}
	\]
	which contradicts Lemma \ref{le2.2}.
	
	If
	\[
	\frac{N+\alpha}{N}<q<p=\frac{N+\alpha+2s}{N},
	\]
	then the vanishing gives
	\[
	Q(u_n)=o_n(1),
	\qquad
	P(u_n)=o_n(1).
	\]
	Therefore,
	\[
	J(u_n)
	=
	\frac12\|u_n\|^2
	+\frac12\int_{\mathbb R^N}V(x)|u_n|^2\,dx
	+o_n(1)\geq o_n(1),
	\]
	which contradicts \(J(u_n)\to m(a)<0\). Thus \eqref{nne2.18} holds.
	
	Since \(u_n\to0\) in \(L^2_{\mathrm{loc}}(\mathbb R^N)\), the sequence \(\{y_n\}\) in \eqref{nne2.18} must satisfy
	\[
	|y_n|\to\infty .
	\]
	Define
	\[
	v_n(x)=u_n(x+y_n).
	\]
	Then \(v_n\in S(a)\), and by the translation invariance of \(J_\infty\),
	\[
	J_\infty(v_n)=J_\infty(u_n)\to m(a).
	\]
	By Lemma \ref{le2.5},
	\[
	m(a)\leq m_\infty(a).
	\]
	On the other hand, by the definition of \(m_\infty(a)\),
	\[
	m_\infty(a)\leq J_\infty(v_n)=m(a)+o_n(1).
	\]
	Hence
	\[
	m_\infty(a)=m(a),
	\qquad
	J_\infty(v_n)\to m_\infty(a).
	\]
	
	By \eqref{nne2.18}, up to a subsequence,
	\[
	v_n\rightharpoonup v \quad \text{in } H^s(\mathbb R^N)
	\]
	for some \(v\not\equiv0\). Applying the argument of Case \(1\) to the autonomous functional \(J_\infty\), we obtain
	\[
	v_n\to v \quad \text{in } H^s(\mathbb R^N),
	\qquad
	v\in S(a),
	\qquad
	J_\infty(v)=m_\infty(a).
	\]
	Finally, since \(V(x)\leq V_\infty\), \(f(x)\geq f_\infty\), and \(g(x)\geq g_\infty\), we have
	\[
	J(v)\leq J_\infty(v).
	\]
	Since \(v\in S(a)\), it follows that
	\[
	m(a)\leq J(v)\leq J_\infty(v)=m_\infty(a)=m(a).
	\]
	Thus
	\[
	J(v)=m(a).
	\]
	Therefore alternative \((II)\) holds. The proof is complete.
\end{proof}

\section{Proof of the main results}

{\bf The proof of Theorem \ref{thm1}:}
By the definition of \(m(a)\), there exists a minimizing sequence
\(\{u_n\}\subset S(a)\) such that
\[
J(u_n)\to m(a)\qquad\text{as }n\to+\infty .
\]
By Lemma \ref{le2.6}, either \(\{u_n\}\) converges strongly in
\(H^s(\mathbb R^N)\), or there exist \(\{y_n\}\subset\mathbb R^N\) with
\(|y_n|\to\infty\) such that \(u_n(\cdot+y_n)\) converges strongly in
\(H^s(\mathbb R^N)\) to some function \(u\in S(a)\) satisfying
\(J(u)=m(a)\). In both cases, there exists \(u\in S(a)\) such that
\[
J(u)=m(a).
\]

Since \(u\) is a constrained minimizer of \(J\) on \(S(a)\), there exists
\(\lambda\in\mathbb R\) such that
\[
J'(u)=\lambda\Psi'(u),
\qquad
\Psi(u)=\frac12\int_{\mathbb R^N}|u|^2\,dx .
\]
Therefore \(u\) is a weak solution of
\[
(-\Delta)^s u+V(x)u
=\lambda u+f(x)\left(I_\alpha *(f|u|^q)\right)|u|^{q-2}u
+g(x)\left(I_\alpha *(g|u|^p)\right)|u|^{p-2}u .
\]

It remains to estimate \(\lambda\). Since
\[
q=\frac{N+\alpha}{N}<p,
\]
we have
\[
\begin{aligned}
	J(u)-\frac{1}{2q}J'(u)u
	=&\left(\frac12-\frac{1}{2q}\right)\|u\|^2
	+\left(\frac12-\frac{1}{2q}\right)
	\int_{\mathbb R^N}V(x)|u|^2\,dx+\left(\frac{1}{2q}-\frac{1}{2p}\right)P(u)>0 .
\end{aligned}
\]
Since \(J'(u)u=\lambda a\), it follows that
\[
m(a)-\frac{\lambda a}{2q}>0.
\]
Hence
\[
\lambda<\frac{2q}{a}m(a).
\]
By Lemma \ref{le2.2},
\[
m(a)<-\frac{f_\infty^2N}{2(N+\alpha)}
S_\alpha^{-\frac{N+\alpha}{N}}a^{\frac{N+\alpha}{N}} .
\]
Consequently,
\[
\lambda
<-f_\infty^2S_\alpha^{-\frac{N+\alpha}{N}}a^{\frac{\alpha}{N}} .
\]
The proof is complete.
\hfill{$\Box$}

{\bf The proof of Theorem \ref{thm2}:}
Let \(\{u_n\}\subset S(a)\) be a minimizing sequence for \(m(a)\). By Lemma
\ref{le2.6}, there exists \(u\in S(a)\) such that
\[
J(u)=m(a).
\]
Thus \(u\) is a constrained minimizer of \(J\) on \(S(a)\). Hence there exists
\(\lambda\in\mathbb R\) such that
\[
J'(u)=\lambda\Psi'(u),
\qquad
\Psi(u)=\frac12\int_{\mathbb R^N}|u|^2\,dx .
\]
Therefore \(u\) is a weak solution of \eqref{eq1}.

Since
\[
\frac{N+\alpha}{N}<q<p=\frac{N+\alpha+2s}{N},
\]
we have
\[
\begin{aligned}
	J(u)-\frac{1}{2q}J'(u)u
	=&\left(\frac12-\frac{1}{2q}\right)\|u\|^2
	+\left(\frac12-\frac{1}{2q}\right)
	\int_{\mathbb R^N}V(x)|u|^2\,dx  \\
	&+\left(\frac{1}{2q}-\frac{1}{2p}\right)P(u)>0 .
\end{aligned}
\]
Since \(J'(u)u=\lambda a\), we obtain
\[
m(a)-\frac{\lambda a}{2q}>0.
\]
Therefore
\[
\lambda<\frac{2q}{a}m(a).
\]
By Lemma \ref{le2.2}, \(m(a)<0\). Hence
\[
\lambda<0.
\]
The proof is complete.
\hfill{$\Box$}

For equation \eqref{neq1}, the associated energy functional is
\[
\begin{aligned}
	I(u)=&
	\frac12\int_{\mathbb R^N}|(-\Delta)^{\frac{s}{2}}u|^2\,dx
	+\frac12\int_{\mathbb R^N}V(x)u^2\,dx
	-\frac{1}{2q}\int_{\mathbb R^N}
	\left(I_\alpha*(f|u|^q)\right)f|u|^q\,dx  \\
	&-\frac{1}{2k}\int_{\mathbb R^N}
	\left(I_\alpha*|u|^k\right)|u|^k\,dx
	-\frac{1}{2p}\int_{\mathbb R^N}
	\left(I_\alpha*(g|u|^p)\right)g|u|^p\,dx .
\end{aligned}
\]
We also set
\[
R(u)=\int_{\mathbb R^N}\left(I_\alpha*|u|^k\right)|u|^k\,dx
\]
and
\[
m_I(a)=\inf_{u\in S(a)}I(u).
\]

{\bf The proof of Theorem \ref{thm3}:}
Since
\[
q=\frac{N+\alpha}{N},\qquad q<k<p=\frac{N+\alpha+2s}{N},
\]
we have
\[
q\gamma_{q,s}=0,\qquad k\gamma_{k,s}<1,\qquad p\gamma_{p,s}=1 .
\]
By \eqref{1.3}, for every \(u\in S(a)\),
\[
\begin{aligned}
	I(u)
	\geq&
	\left(\frac12-\frac{g_{\max}^2C_{\alpha,p}}{2p}a^{p-1}\right)\|u\|^2
	-\frac{f_{\max}^2}{2q}S_\alpha^{-q}a^q  \\
	&-\frac{C_{\alpha,k}}{2k}a^{k(1-\gamma_{k,s})}\|u\|^{2k\gamma_{k,s}} .
\end{aligned}
\]
Since
\[
0<a<a^*=\left(\frac{p}{C_{\alpha,p}g_{\max}^2}\right)^{\frac{1}{p-1}},
\]
the coefficient of \(\|u\|^2\) is positive. Moreover,
\(2k\gamma_{k,s}<2\). Thus \(I\) is bounded from below and coercive on
\(S(a)\).

Next we prove that \(m_I(a)\) is strictly below the lower critical level.
Let \(v\) be an extremal function for \(S_\alpha\), and set
\[
u=\frac{\sqrt a\,v}{\|v\|_2},
\qquad
u_t(x)=t^{\frac N2}u(tx),\quad t>0 .
\]
Then \(u_t\in S(a)\). Moreover,
\[
\int_{\mathbb R^N}
\left(I_\alpha*|u|^q\right)|u|^q\,dx
=a^qS_\alpha^{-q}.
\]
A direct computation gives
\[
\begin{aligned}
	I(u_t)
	\leq&
	\frac12t^{2s}\|u\|^2+\frac{\|V\|_\infty}{2}a
	-\frac{f_\infty^2}{2q}a^qS_\alpha^{-q}-\frac{1}{2k}t^{Nk-N-\alpha}R(u)
	-\frac{g_\infty^2}{2p}t^{Np-N-\alpha}\widetilde P(u),
\end{aligned}
\]
where
\[
\widetilde P(u)=
\int_{\mathbb R^N}
\left(I_\alpha*|u|^p\right)|u|^p\,dx .
\]
Since \(q<k<p\), we have
\[
0<Nk-N-\alpha<2s.
\]
Hence, for \(t>0\) sufficiently small,
\[
\frac12t^{2s}\|u\|^2
-\frac{1}{2k}t^{Nk-N-\alpha}R(u)<0.
\]
Choosing \(V_*=V_*(a)>0\) sufficiently small, we obtain, for every
\(\|V\|_\infty<V_*\),
\[
\frac12t^{2s}\|u\|^2
-\frac{1}{2k}t^{Nk-N-\alpha}R(u)
+\frac{\|V\|_\infty}{2}a<0.
\]
Therefore
\[
m_I(a)\leq I(u_t)
<-\frac{f_\infty^2}{2q}a^qS_\alpha^{-q}
=-\frac{f_\infty^2N}{2(N+\alpha)}
S_\alpha^{-\frac{N+\alpha}{N}}a^{\frac{N+\alpha}{N}} .
\]

Let \(\{u_n\}\subset S(a)\) be a minimizing sequence for \(m_I(a)\). By the
coercivity of \(I\), \(\{u_n\}\) is bounded in \(H^s(\mathbb R^N)\). The compactness
argument in Lemma \ref{le2.6} applies to \(I\) with only one additional term,
namely \(R\). This term satisfies the nonlocal Brezis--Lieb splitting
\[
R(u_n)=R(u_n-u)+R(u)+o_n(1)
\]
whenever \(u_n\rightharpoonup u\) in \(H^s(\mathbb R^N)\) and \(u_n\to u\)
a.e. in \(\mathbb R^N\). It is also invariant under translations. Hence the
dichotomy exclusion and the compactness alternative in Lemma \ref{le2.6}
remain valid for \(I\).

More precisely, the corresponding limiting functional is
\[
\begin{aligned}
	I_\infty(u)=&
	\frac12\|u\|^2
	+\frac{V_\infty}{2}\int_{\mathbb R^N}u^2\,dx
	-\frac{f_\infty^2}{2q}\int_{\mathbb R^N}
	\left(I_\alpha*|u|^q\right)|u|^q\,dx  \\
	&-\frac{1}{2k}\int_{\mathbb R^N}
	\left(I_\alpha*|u|^k\right)|u|^k\,dx
	-\frac{g_\infty^2}{2p}\int_{\mathbb R^N}
	\left(I_\alpha*|u|^p\right)|u|^p\,dx ,
\end{aligned}
\]
and
\[
m_{I,\infty}(a)=\inf_{u\in S(a)}I_\infty(u).
\]
Since \(V(x)\leq V_\infty\), \(f(x)\geq f_\infty\), and \(g(x)\geq g_\infty\), we have
\[
I(u)\leq I_\infty(u),\qquad u\in S(a).
\]
Thus, as in Lemma \ref{le2.6}, the possible loss of compactness at infinity
leads to a translated minimizing sequence for \(m_{I,\infty}(a)\). In this case,
the comparison \(I(u)\leq I_\infty(u)\) on \(S(a)\) gives
\[
m_I(a)\leq m_{I,\infty}(a),
\]
while the translated minimizing sequence gives
\[
m_{I,\infty}(a)\leq m_I(a).
\]
Hence
\[
m_{I,\infty}(a)=m_I(a).
\]
Therefore there exists \(v\in S(a)\) such that
\[
I(v)\leq I_\infty(v)=m_{I,\infty}(a)=m_I(a).
\]
Since \(m_I(a)\leq I(v)\), we obtain
\[
I(v)=m_I(a).
\]
Consequently, \(m_I(a)\) is achieved by some \(u\in S(a)\).

Since \(u\) is a constrained minimizer of \(I\) on \(S(a)\), there exists
\(\lambda\in\mathbb R\) such that
\[
I'(u)=\lambda\Psi'(u),
\qquad
\Psi(u)=\frac12\int_{\mathbb R^N}|u|^2\,dx .
\]
Therefore \(u\) is a weak solution of
\[
\begin{aligned}
	(-\Delta)^s u+V(x)u
	=&\lambda u
	+f(x)\left(I_\alpha*(f|u|^q)\right)|u|^{q-2}u+\left(I_\alpha*|u|^k\right)|u|^{k-2}u
	+g(x)\left(I_\alpha*(g|u|^p)\right)|u|^{p-2}u .
\end{aligned}
\]

Finally, since \(q<k<p\), we have
\[
\begin{aligned}
	I(u)-\frac{1}{2q}I'(u)u
	=&\left(\frac12-\frac{1}{2q}\right)\|u\|^2
	+\left(\frac12-\frac{1}{2q}\right)
	\int_{\mathbb R^N}V(x)|u|^2\,dx  \\
	&+\left(\frac{1}{2q}-\frac{1}{2k}\right)R(u)
	+\left(\frac{1}{2q}-\frac{1}{2p}\right)P(u)>0 .
\end{aligned}
\]
Since \(I'(u)u=\lambda a\), we obtain
\[
m_I(a)-\frac{\lambda a}{2q}>0.
\]
Therefore
\[
\lambda<\frac{2q}{a}m_I(a).
\]
Using the estimate for \(m_I(a)\), we derive
\[
\lambda
<-f_\infty^2S_\alpha^{-q}a^{q-1}
=-f_\infty^2S_\alpha^{-\frac{N+\alpha}{N}}a^{\frac{\alpha}{N}} .
\]
The proof is complete.
\hfill{$\Box$}

\section*{Acknowledgment}
We express our gratitude to the anonymous referee for their meticulous review of our manuscript and valuable feedback provided for its enhancement. 

Y. Chen was supported by the National Natural Science Foundation of China (12161007) and Guangxi Natural Science Foundation Project (2023GXNSFAA026190).
Z. Yang was supported by the National Natural Science Foundation of China (12301145, 12261107), Yunnan Fundamental Research Projects (202201AU070031, 202401AU070123) and Yunnan Revitalization Talent Support Program. J. Zhang was supported by National Natural Science Foundation of China (12371109).

\end{document}